\documentclass[11pt,leqno]{article}

\usepackage{amsthm}
\usepackage{amsmath}
\usepackage{amssymb}
\usepackage{amsfonts}
\usepackage{amscd}
\usepackage{palatino}
\usepackage[mathscr]{eucal}
\usepackage{epsfig}
\usepackage{verbatim}
\usepackage{latexsym}
\usepackage{graphics}

\begin{document}

\newcounter{chn}
\newenvironment{nch}[1]
               {\addtocounter{chn}{1}
               \newpage
                      \begin{center}
                      {\Large {\bf {#1}}}
                      \end{center}
                      \vspace{5mm}}{}

%The Theorem env.
\newcounter{thn}
\renewenvironment{th}
{\addtocounter{thn}{1}
 \em 
     \noindent 
     {\sc Theorem \arabic{thn} . ---}}{}

%Theorem without numbers
\newenvironment{theo}
   {\em
        \vspace{2mm}
        \noindent
                {\sc Theorem\/. ---}}{}

%Lemma without numbers
\newenvironment{lem}
   {\em
        \vspace{2mm}
        \noindent
                {\sc Lemma\/. ---}}{}

%Prop without numbers
\newenvironment{pro}
   {\em
        \vspace{2mm}
        \noindent
                {\sc Proposition\/. ---}}{}

%The Lemma env.
\newcounter{len}
\renewenvironment{le}
{\addtocounter{len}{1}
 \em 
     \vspace{2mm}
     \noindent 
     {\sc Lemma \arabic{len}. ---}}{}

%The Proposition env.
\newcounter{prn}
\newenvironment{prop}
{\addtocounter{prn}{1}
  \em 
      \vspace{2mm}
      \noindent 
      {\sc Proposition \arabic{prn}. ---}}{}

\newcounter{exn}
\newenvironment{prob}
                     {\addtocounter{exn}{1} 
                      \vspace{2mm}
                      \noindent 
                      {\em Exercise \arabic{exn} :}}{}

\newcounter{corn}
\newenvironment{cor}
  {\em 
     \addtocounter{corn}{1}
     \vspace{2mm}   
     \noindent 
        {\sc Corollary \arabic{corn}. ---}}{\vspace{1mm}}

%Corollary without numbers
\newenvironment{cor1}
   {\em 
        \vspace{2mm}
        \noindent
                {\sc Corollary\/. ---}}{}

%Definitions
\newcounter{defn}
\newenvironment{df}{
                    \vspace{2mm}
                    \addtocounter{defn}{1} 
                    \noindent 
                            {\sc Definition \arabic{defn}. ---}}{\vspace{1mm} }{}
%paragraphs
\newcounter{pgn}
\newenvironment{pg}{
                    \vspace{2mm}
                    \addtocounter{pgn}{1}
                    \noindent
                            {\arabic{scn}.\arabic{pgn}\;}}{}

%The section env.
\newcounter{scn}
\newenvironment{nsc}[1]
              {\addtocounter{scn}{1}
              \vspace{3mm}
                     \begin{center}
                     {\sc { \arabic{scn}. {#1}}}
                     \end{center}
              \setcounter{equation}{0}
              \setcounter{prn}{0}
              \setcounter{len}{0}
              \setcounter{thn}{0}
              \setcounter{pgn}{0}
              \setcounter{corn}{0}
              \vspace{3mm}}{}

\begin{center}
\Large {\bf

{
The  Large Sieve Inequality for Integer Polynomial Amplitudes
}
}

\end{center}
\vspace{5mm}

\vspace{-10mm}
\begin{center}
\it{

by \\

\vspace{3mm}

Gyan Prakash and D.S. Ramana

}
\end{center}

\vspace{2mm}
\begin{abstract}
\noindent We obtain a close to the best possible version of the large sieve inequality with amplitudes given by the values of a polynomial with integer coefficients of degree $\geq 2$.
\end{abstract}

\begin{nsc}{Introduction}

\noindent
It is of interest in the context of inequalities of the large sieve type to obtain estimates for the sum $\sum_{x \in \mathcal{X}}  |\sum_{i \in I } a_i e(xf(i))|^2$, where $e(z)$ denotes $e^{2 \pi i z}$ for any complex number $z$, $\mathcal{X}$ is a well-spaced set of real numbers, $I$ is an interval of the integers , $\{a_i\}_{i \in I}$ are complex numbers
and $f$ is real valued function on $I$ such that $f(I)$ is sparse, that is, the length of $f(I)$ is much larger than the length of $I$. When $f(I)$ is sparse , the duality argument that is used to establish the classical large sieve inequality generally gives weak bounds thus provoking the search for alternate arguments. 

\vspace{2mm}
\noindent
Basic example of sparse sequences are the sequence of values of polynomials of degree $\geq 2$ and Iwaniec and Kowalski, in their book \cite{iwko} (see page 184, the paragraph following Problem 7.19), pose the problem of determining good large sieve bounds when $f$ is $P(T)$, a given polynomial in ${\bf Z}[T]$. Our purpose in this note is to verify the following theorem, which provides a result in this direction.

\vspace{2mm}
\noindent
In the theorem below and thereafter $\omega(n)$ denotes the number of prime divisors of $n$ and  $\|a\|^2$ denotes $\sum_{i \in I} |a_i|^2$ for a finite sequence of complex numbers $\{a_i\}_{i \in I}$. Further, for each integer $k\geq 1$, we define

\vspace{-3mm}
\noindent
\begin{equation}
\label{tk}
\theta(k) =  k\binom{k+1}{2}  \; .
\end{equation}

\vspace{2mm}
\begin{th} 
Let $Q$ and $k$ be integers $\geq 1$ and $I$ be an interval in ${\bf Z}$ of length $N$. When ${\mathcal F}(Q)$ is 
the Farey sequence of order $Q$ and $P(T) = c_0T^k + c_1 T^{k-1} + \ldots + c_k$ is a polynomial of degree $k$ in ${\bf Z}[T]$ we have the inequality 
                                                                                    \begin{equation}  
\label{lsth1}
\sum_{x \in {\mathcal F}(Q)}
\left| \sum_{ i \in I } a_i e( x P(i))\right|^2
\; \ll \; 
Q(N + Q) (\log Q)^{\omega(c_0) + \theta(k)} \, \|a\|^2 \;  
\end{equation}

\noindent
for every sequence of complex numbers $\{a_i\}_{i \in I}$, where the constant implicit in the $\ll$ depends only on $k$.  

\end{th}

\vspace{2mm}
\noindent
While the classical large sieve inequality gives a much better bound than that given by \eqref{lsth1} when $P(T)$ is linear, the bound \eqref{lsth1} is the best possible upto the term $(\log Q)^{\omega(c_0) + \theta(k)}$ and constant implicit in the $\ll$ in \eqref{lsth1} when the degree of $P(T)$ is $\geq 2$. We show this by means of an example in Section 3, where we provide a proof of Theorem 1.

\vspace{2mm}
\noindent
A number of authors (see \cite{ram}, \cite{zh}, \cite{st}) have recently obtained upper bounds for the sum on the left hand side of \eqref{lsth1} from various points of view. These bounds are, however, are comparable to that  given by \eqref{lsth1} only when $P(T)$ is of degree 2 and the interval $I$ is of the form $(M, M+N]$ with $M \ll N$. In fact, the only bound for polynomials $P(T)$ of degree $\geq 3$ that we are aware of is due to S. Baier who uses the method of Zhao \cite{zh} method to observe that (see Corollary 3 following Theorem 2 in \cite{st}~) when $P(T)= T^m$, for any integer $m \geq 3$, and when $I$ is of the form $(0, N]$ then the left hand side of \eqref{lsth1} is $\ll_{\epsilon} (NQ^{2(1 -1/m)} + Q^2) N^{1+\epsilon} \sup_{0 < i \leq N} |a_i|^2 $ under Hooley's hypothesis $K^{*}$ in the context of Waring's problem. Baier deduces this from an estimate that is valid even when the Farey series is replaced with an arbitrary well spaced sets as well. 

\vspace{2mm}
\noindent
When indeed $P(T)$ is of degree 2 and the interval $I$ is of the form $(0, N]$, Ramar{\'e}`s method, based on his theory of local factors, gives the bound $Q(N+Qg(Q))(\log_2 2Q)^2$ for the sum on the left hand side \eqref{lsth1}, where $g(Q) = \exp(C\log_2 Q \log_3 Q)$. On the other hand, Zhao\cite{zh} gives,  for the same sum, upper bounds essentially of the form $(Q(NM)^{1/2} + Q^2)(NM)^{\epsilon}$, for each $\epsilon >0$ when $I$ is of the form $(M,M+N]$, via an elegant application of the double large sieve inequality. While Ramar{\'e}'s estimate is sharper than that given by Theorem 1 when $N$ is much larger than $Q$, Zhao's estimate has the advantage of being applicable even when the Farey series is replaced with an arbitrary well spaced set.

\vspace{2mm}
\noindent
In contrast to the aforementioned results, Theorem 1 is valid for {\em all} integer polynomials $P(T)$ and the bound given by this theorem is {\em uniform} with regard to the position of the interval $I$.

\end{nsc}

\begin{nsc}{Number of Zeros of $P(T)$ modulo $m$}

\noindent
Let $P(T) = c_0T^k + c_1T^{k-1} + \ldots + c_{k}$ be a polynomial in ${\bf Z}[T]$ and, for any integer $m \geq 1$, let $S(m)$ be the set of congruence classes $l$ modulo $m$ such that $P(l) \equiv 0 \,{\rm mod}\, m$ and let $\rho(m)$ be   ${\rm Card}(S(m))$. Let $Q$ be a real number $\geq 1$. Proposition 1 below gives for $\sum_{1 \leq m \leq Q} \frac{\rho(m)}{m}$ an upper bound that is independent of the constant term of $P(T)$. The proof of Theorem 1 relies crucially on this feature of this bound.

\vspace{2mm}
\noindent
Let $p$ be a prime number and $m$, $n$ be integers $\geq 1$. When $ m \geq n$ the image of $S(m)$ under the canonical map from ${\bf Z}/p^m{\bf Z}$ onto ${\bf Z}/p^n{\bf Z}$ is contained in $S(n)$. Therefore we have $\frac{\rho(p^m)}{p^m} \leq \frac{\rho(p^n)}{p^n}$ whenever $m \geq n \geq 1$.

\vspace{2mm}
\noindent
Suppose now that $p$ is a prime number that does not divide $c_0$. We then have $\rho(p) \leq k$ and hence $\frac{\rho(p^m)}{p^m} \leq \frac{k}{p}$ for all $m \geq 1$. We shall presently improve upon this upper bound for large $m$. To this end we set, for any integer $m \geq 1$,  $a(m,k)$ to denote the smallest integer $\geq \frac{m}{\binom{k+1}{2}}$ and we verify that any interval of the real line of length $p^{a(m,k)}$ contains no more than $k+2$ integers $x$ such that $P(x)$ is divisible by $p^m$. To verify this it suffices to show that when $x_1 , x_2, \ldots, x_{k+1}$ are distinct integers such that $P(x_i)$ is divisible by 
$p^m$ for each $i$, we have $\sup_{(i,j)} |x_i - x_j| \geq p^{ a(m,k)}$. Indeed, on recalling the well known identity for the vandermonde determinant we have

\vspace{-3mm}
\begin{equation}
\label{vander}
c_0\prod_{1 \leq i < j \leq k+1} (x_i - x_j)
\;=\;
\left| \begin{matrix}
1 & 1 & \dots & 1\\
{x_1} & {x_2} & \dots & {x_{k+1}}\\
\vdots & \vdots & & \vdots \\
{x_1^{k-1}} & {x_2^{k-1}} & \dots & {x_{k+1}^{k-1}}\\
P(x_1) & P(x_2) & \dots & P(x_{k+1})\\
\end{matrix} \right | \; .
\end{equation} 

\noindent
Since the right hand side of \eqref{vander} is divisible by $p^m$ and $p$ does not divide $c_0$ we see that $p^m$ divides $\prod_{1 \leq i < j \leq k+1} (x_i - x_j)$. Consequently, 

\vspace{-3mm}
\begin{equation}
\label{vander1}
{\binom{k+1}{2}} \sup_{i \neq j} v_p(x_i - x_j) \geq \sum_{1 \leq i < j \leq k+1} v_p(x_i - x_j) \geq m,
\end{equation}

\noindent
from which it follows that $\sup_{i\neq j} v_p(x_i - x_j) \geq {a(m,k)}$ and, because the $x_i$ are distinct, that $\sup_{i,j}|x_i - x_j| \geq p^{a(m,k)}$.   

\vspace{2mm}
\noindent
For each integer $m \geq 1$, the set $S(m)$ is in bijection with the subset of
the integers $x$ in the interval $[0, p^m)$ such that $P(x)$ is divisible by $p
^m$. On dividing this interval into subintervals of length $ p^{a(m,k)}$ we then conclude that when $p$ does not divide $c_0$ we have $\rho(p^m) \leq \frac{(k+2) p^m}{p^{a(m,k)}}$, for all $m \geq 1$.     

\vspace{2mm}
\noindent
With the aid of the bounds for $\frac{\rho(p^m)}{p^m}$ given above we then conclude that when $p$ is a prime number that does not divide $c_0$ we have 

\vspace{-3mm}
\begin{equation}
\label{bound1}
\sum_{m \geq 0} \frac{\rho(p^m)}{p^m} \; \leq \; 1 + \sum_{1 \leq m \leq \binom{k+1}{2}} \frac{k}{p} + \sum_{ m > \binom{k+1}{2}} \frac{k+2}{p^{a(m,k)}} \; = \; 1 + \frac{\theta(k)}{p} \; + \; (k+2)\binom{k+1}{2}\sum_{m \geq 2}\frac{1}{p^m},
\end{equation} 

\noindent
where the last identity follows on dividing the sum over $m \geq \binom{k+1}{2}$ into sums over congruence classes modulo $\binom{k+1}{2}$ and noting that $a(l+d\binom{k+1}{2},k) = d$, when $d$ is any integer and $l$ an integer satisfying $0\leq l < \binom{k+1}{2}$.

\begin{prop} Let $P(T) = c_0T^k + c_1T^{k-1} + \ldots + c_{k}$ be a polynomial of degree $k$ in ${\bf Z}[T]$ and let $\rho(m)$, for each integer $m \geq 1$, be the number of residue classes $l$ modulo $m$  such that $P(l) \equiv \, 0 \,{\rm mod}\, m$. For any real number $Q \geq 1$ we then have $\sum_{1 \leq m \leq Q} \frac{\rho(m)}{m} \ll (\log Q)^{\omega(c_0) +  \theta(k)}$, where the implied constant in the $\ll$ is depends only on $k$.
\end{prop}

\vspace{2mm}
\noindent
{\sc Proof. ---} Since $\rho(m)$ is multiplicative, the sum  $\sum_{1 \leq m \leq Q} \frac{\rho(m)}{m}$ is evidently majorised by  

\vspace{-3mm}
\begin{equation}
\label{bound2}
\prod_{1 \leq p \leq Q}  \left(\sum_{m \geq 0, \atop p^m \leq Q} \frac{\rho(p^m)}{p^m} \right) \; \leq \; (\log Q)^{w(c_0)} \prod_{1 \leq p \leq Q, \atop (c_0,p) =1 } \left( 1 + \frac{\theta(k)}{p} + (k+2) \binom{k+1}{2}\sum_{m\geq 2}\frac{1}{p^m} \right) , 
\end{equation}

\noindent
 where we have used \eqref{bound1} when $p$ does not divide $c_0$ and $\rho(p^m) \leq p^m$ otherwise. The proposition now follows on dropping the condition $(c_0,p) =1$ in the product on the right hand side of \eqref{bound2} and noting that $\prod_{1 \leq p \leq Q} (1 + \frac{a}{p} + \frac{b}{p^2}) \ll_b (\log Q)^a$, for any real numbers $a$ and $b$.
\end{nsc}

\begin{nsc}{Proof of the Large Sieve Inequality}

\noindent
For each $(i,j) \in I \times I $ let us set 

\vspace{-2mm}
\begin{equation}
\label{eq0.5}
K(i,j) \, = \, \sum_{x \in {\mathcal F}(Q)} e( x (P(i)-P(j)) \; 
\end{equation}

\noindent
so that $|K(i,j)| = |K(j,i)|$, for each $(i,j)$. On squaring out the sum over $i \in I$ in the left hand side of the inequality given by Theorem 1, interchanging the summations and applying the triangle inequality together with $|a_i \bar{a_j}| \leq \frac{1}{2}(|a_i|^2 + |a_j|^2)$ for each $(i,j)$ we obtain

\vspace{-2mm}
\begin{equation}
\label{eq1}
\sum_{x \in {\mathcal F}(Q)}
\left| \sum_{ i \in I } a_i e( x P(i))\right|^2
= 
\sum_{ (i,j) \in I \times I } a_i \bar{a_j} K(i,j) \; \leq \; (\sup_{j \in I} \sum_{ i \in I} |K(i,j)|)\, \|a\|^2 \; .
\end{equation} 

\noindent
Let $c(i,j)$ denote $P(i) -P(j)$ for each $(i,j)$. The classical estimate $|\sum_{0 \leq p \leq q-1, \atop (p,q) =1 .} e(\frac{ap}{q})| \, \leq \, (a,q)$, valid for any integer $a$ with the convention that $(0,q) =q$, then implies  

\vspace{-2mm}
\begin{equation}
\label{eq3}
|K(i,j)| \; \leq\; \sum_{1 \leq q \leq Q} |\sum_{0 \leq p \leq q-1, \atop (p,q) =1 .} e(\frac{pc(i,j)}{q})| \; \leq\; \sum_{1 \leq q \leq Q} (c(i,j), q) \; . 
\end{equation}

\noindent
Since for any integer $k$ with $1 \leq k \leq Q$, the number of multiples $q$ of $k$ with $1 \leq q \leq Q$ does not exceed $\frac{Q}{k} + 1 \leq \frac{2Q}{k}$, we obtain  

\vspace{-2mm}
\begin{equation}
\label{eq4}
\sum_{1 \leq q \leq Q} (c(i,j), q) \; \leq \;
\sum_{1 \leq k \leq Q, \atop k | c(i,j)} k  \sum_{ 1 \leq q \leq Q, \atop q \equiv 0 {\rm {mod}} k } 1 \; \leq \; 2Q\sum_{1 \leq k \leq Q, \atop k | c(i,j)} 1.
\end{equation}

\noindent
For any $j \in I$, let us set $\rho_{j}(k)$ to denote the number of congruence classes $l$ modulo $k$ for which $P(l) \equiv P(j) \, {\rm mod}\, k$.
On combining \eqref{eq4} with \eqref{eq3} and recalling that $I$ is an interval of length $N$, we then conclude that for each $j \in I$  

\vspace{-2mm}
\begin{equation}
\label{eq5}
\sum_{ i \in I} |K(i,j)| \, \leq \, 2Q\sum_{ i \in I}  \sum_{1 \leq k \leq Q, \atop k | c(i,j)} 1 = 2Q \sum_{1 \leq k \leq Q} \sum_{ i \in I, \atop {c(i,j) \equiv 0 \,{\rm mod}\, k} } 1 \leq  2Q \sum_{1 \leq k \leq Q}  \rho_{j}(k) \left(\frac{N}{k} + 1\right) \; .
\end{equation}  

\noindent
On applying Proposition 1 of Section 2 to the polynomial $P(T) - P(j)$, we see, for each $j \in I$, that

\vspace{-2mm}
\begin{equation}
\label{eq6}
 \sum_{1 \leq k \leq Q}  \rho_{j}(k) \left(\frac{N}{k} + 1\right) \; \leq \; \sum_{1 \leq k \leq Q} (N+Q) \frac{\rho_{j}(k)}{k} \ll (N+Q)(\log Q)^{\omega(c_0 ) + \theta(k)} \; ,
\end{equation}

\noindent
combining which with \eqref{eq1} and \eqref{eq5} we obtain Theorem 1.

\vspace{2mm}
\noindent
Let us verify that upto the term $(\log Q)^{\omega(c_0) + \theta(k)}$ and the constant implicit in $\ll$ the bound given by Theorem 1 is the best possible. To this end, let us take $P(T) = T^n$, where $n \geq 2$. We then learn on page 24 of \cite{kor} that when $q$ is a {\em prime number} we have

\vspace{-3mm}
\begin{equation}
\label{ex1}
\sum_{1\leq p \leq q-1}\left| \sum_{1 \leq i \leq q} e(\frac{pP(i)}{q}) \right|
^2 \; = \; (n-1)q(q-1) \; .
\end{equation}

\noindent
Moreover, we  have the bound $|\sum_{1 \leq i \leq q} e(\frac{pP(i) + ki}{q})|
\leq (n-1)q^{1/2}$ from the estimate of Weil for exponential sums, for all prime numbers $q > n$ and all integers $p, k$ with $(p,q) = 1$. On using Theorem 2,
 page 12 of \cite{kor} we then deduce that

\vspace{-3mm}
\begin{equation}
\label{ex1.5}
\left|\sum_{1 \leq i \leq m} e(\frac{pP(i)}{q})\right| \leq 2(n-1)q^{1/2}\log q
  ,
\end{equation}

\noindent
for all prime numbers $q > n$ and all integers $p, m$ with $(p,q) = 1$ and $1 \leq m \leq q$. Let us now take $N$ and $Q$ integers $\geq 1$ with $N \geq 8(n-1
)Q\log Q$ and $Q \geq n^2$. On dividing the interval $(0, N]$ into subintervals
 of length $q$ and applying the triangle inequality together with \eqref{ex1} we then see that

\vspace{-3mm}
\begin{equation}
\label{ex2}
\sum_{1\leq p \leq q-1}
\left| \sum_{1 \leq i \leq N} e(\frac{pP(i)}{q}) \right|^2  \geq
\sum_{1\leq p \leq q-1} \, \left(\left[\frac{N}{q}\right] \left| \sum_{1 \leq i
 \leq q} e(\frac{pP(i)}{q}) \right| - 2(n-1)q^{1/2}\log{ q }\right)^2  \geq    \frac{(n-1)N^2}{8} ,
\end{equation}

\noindent
for every {\em prime number} $q$ satisfying $n < q \leq Q$. Consequently, we ob
tain the minorisation

\vspace{-3mm}
\begin{equation}
\label{ex3}
\sum_{x \in {\mathfrak F}(Q)} \left| \sum_{ 0 < i \leq N }  e( xP(i))\right|^2
\; \geq \;
\sum_{1 \leq q \leq Q, \atop {q  {\; \text prime}, \atop q > n .}} \sum_{1 \leq p \leq q-1} \left| \sum_{ 0 < i \leq N }  e( \frac{pP(i)}{q})\right|^2 \; \geq \;
\frac{(n-1)N^2Q}{16\log Q} \; ,
\end{equation}

\noindent
which may be compared with the upper bound $\ll_{n} N^{2}Q(\log Q)^{\theta(n)}$ for the first term on left hand side of \eqref{ex3} supplied by Theorem 1 when applied with $I$ taken to be the interval of integers $(0,N]$, the polynomial $P(T) = T^n$ and all the $a_i =1$. When $N < 8(n-1)Q\log Q$, this theorem gives the upper bound $\ll_{n} Q^{2}(\log Q)^{\theta(n)+1}\|a\|^2$ for the left hand side of \eqref{lsth1}, which may be compared with the lower bound $Q^2\|a \|$ obtained when $I$ taken to be the interval $(0,1]$ in this expression.

\vspace{2mm}
\noindent
The method of proof of Gallagher's inequality (see page 144 of \cite{ramm}) immediately implies the following corollary to Theorem 1.

\begin{cor1} 
Let $D$ be an integer $\geq 1$. When $I$ is an interval of the integers of length $N$ and $P(T) = c_0T^k + c_1 T^{k-1} + \ldots + c_k$ is a polynomial of degree $k$ in ${\bf Z}[T]$ we have the inequality

\begin{equation}
\label{lschi}
\sum_{1 \leq d \leq D} \frac{\phi(d)}{d}\sum_{ \chi \,{\rm mod}^{*}\, d}
\left| \sum_{ i \in I } a_i \chi(P(i))\right|^2
\; \ll \;
D(N + D) (\log D)^{\omega(c_0) + \theta(k)} \, \|a\|^2 \;
\end{equation}

\noindent
for every sequence of complex numbers $\{a_i\}_{i \in I}$, where the constant  implicit in the $\ll$ depends only on $k$.

\end{cor1}

\end{nsc}

%\begin{nsc}{A Bombieri-Vinogradov Type Theorem}

%\input bv

%\end{nsc}

\vspace{3mm}
\noindent
{\bf Acknowledgement :} Our deepest thanks go to Professor Olivier Ramar{\'e} for generously providing us with his time and suggestions. We are indebted to Dr. Liangyi Zhao, Professor R. Balasubramanian and   Professor R. Heath-Brown for discussing the problem addressed here with us.

\vspace{1cm}
\begin{flushleft}

{\em
Harish-Chandra Research Institute, \\
Chhatnag Road, Jhunsi,\\
Allahabad - 211 019, India.}\\

{\em email} : gyan@mri.ernet.in, suri@mri.ernet.in
\end{flushleft}

\end{document}